\documentclass[12pt,oneside,english]{amsart}
\usepackage{bbding}
\usepackage{dsfont}
\usepackage{graphicx}
\usepackage{amsthm,amscd}
\usepackage{calc}

\usepackage{mathrsfs}
\usepackage{mathrsfs}
\usepackage{CJK,fancyhdr,amscd}
\usepackage{anysize}
\usepackage{amsmath,amssymb,amsfonts}
\usepackage{epsfig}
\usepackage{latexsym}
\usepackage[pagebackref,colorlinks]{hyperref}
\usepackage{indentfirst, latexsym}
\usepackage{graphics}
\usepackage{amsmath, bm}
\usepackage{amsfonts}
\usepackage{amssymb}%
\usepackage{color}
\usepackage{mathtools}
\usepackage{verbatim}

\providecommand{\U}[1]{\protect\rule{.1in}{.1in}}

\numberwithin{equation}{section}
\providecommand{\U}[1]{\protect\rule{.1in}{.1in}}

\newtheorem{theorem} {Theorem} [section]
\newtheorem{proposition}[theorem]{Proposition}
\newtheorem{corollary}  [theorem]     {Corollary}
\newtheorem{lemma}  [theorem]     {Lemma}
\newtheorem{example}  [theorem]     {Example}
\newtheorem{remark}  [theorem]     {Remark}
\newtheorem{definition}  [theorem]     {Definition}

\newcommand{\bproof}{\begin{proof}}
\newcommand{\eproof}{\end{proof}}

\newcommand{\p}{\partial}

\newcommand{\pj}{\partial_J}
\newcommand{\pp}{\partial\partial_J}
\newcommand{\supp}{\mathrm{supp}\,}

\newcommand{\bthm}{\begin{theorem}}
\newcommand{\ethm}{\end{theorem}}
\newcommand{\blem}{\begin{lemma}}
\newcommand{\elem}{\end{lemma}}
\newcommand{\bprop}{\begin{proposition}}
\newcommand{\eprop}{\end{proposition}}
\newcommand{\brmk}{\begin{remark}}
\newcommand{\ermk}{\end{remark}}
\newcommand{\bcor}{\begin{corollary}}
\newcommand{\ecor}{\end{corollary}}
\newcommand{\bdefi}{\begin{definition}}
\newcommand{\edefi}{\end{definition}}
\newcommand{\beq}{\begin{equation}}
\newcommand{\eeq}{\end{equation}}
\newcommand{\bex}{\begin{example}}
\newcommand{\eex}{\end{example}}

\newcommand{\mpsh}{\mbox{QPSH}(\Omega)}
\newcommand{\mpshn}{\mbox{QPSH}^{-}(\Omega)}
\usepackage[lmargin=1in, rmargin=1in]{geometry}
\title{degenerate quaternionic Monge-Amp\`ere equations in weighted energy classes}
\author{Genglong Lin}
\address{Bejing Institute of Mathematical Sciences and Application, Beijing 101408, China}
\email{lingenglong@amss.ac.cn}

\begin{document}
\subjclass[1991]{32U15, 31C10, 32U40}
\date{March 3, 2025}
\keywords{Monge-Amp\`ere equation, quaternionic plurisubharmonic function, Cegrell type classes}
\begin{abstract}
In this paper, we consider degenerate quaternionic Monge-Amp\`ere equations in weighted energy
class $\mathcal{E}_{\chi}(\Omega)$ where $\Omega$ is a quarternionic
domain in $\mathbb{H}^n$ and $\chi$ is a weight function which satisfies some natural conditions. Firstly we prove that the quaternionic Monge-Amp\`ere
operator is well-defined for functions  in $\mathcal{E}_{\chi}(\Omega)$, in particular $\mathcal{E}_p(\Omega),p>0$.

Secondly, we prove that fine property holds in the Cegrell type class $\mathcal{E}(\Omega)$. As an application,
we prove a mass concentration theorem for the quarternionic plurisubharmonic envelope.

In the study of complex Monge-Amp\`ere equation, characterization of finite energy range of complex Monge-Amp\`ere
operator was a central problem which aroused the interest of experts in the subject.
As a quaternionic analogue, we prove a theorem which explicitly characterizes 
the  finite energy range of \emph{quaternionic}  Monge-Amp\`ere operator in the end.
\end{abstract}
\maketitle

\section{introduction}\label{introduction}
Pluripotential theory, initiated in the seminal papers \cite{BT82}\cite{BT79}, plays 
an important role in solving problems of complex analysis and geometry. It has been
generalized in many directions in the last decade, among which is the calibrated 
geometry studied by Harvey-Lawson \cite{HL09}. An interesting direction is the basics
of pluripotential theory in $\mathbb{H}^n$ which established by Alesker \cite{Ale03b}, 
and more generally on hypercomplex manifold  in \cite{AV06}.

Inspired by \cite{BT76},
Alesker developed the foundations of pluripotential theory in the quaternionic setting
and showed that the definition of quaternionic Monge-Amp\`ere operator defined for smooth function as the 
Moore determinant of a quaternionic Hessian can be extended to the class of continuous
functions \cite{Ale03b}. The Dirichlet problem in a quaternionic strictly pseudoconvex
domain $\Omega\subset\mathbb{H}^n$ with a continuous data was solved in \cite{Ale03a}.
Wan \cite{Wan20b} showed the existence of continuous solution to the Dirichlet
problem provided the right-hand-side has $L^p,p\geq4$ densities. This result
was improved into the  case $p\geq2$ in a paper of Sroka \cite{Sro20}.

As observed by Sroka \cite{Sro20}\cite{Sro21}, chosen a suitable frame of $\mathbb{H}^n$ the operator $\partial$ and the twisted
operator $\pj$ coincided with $d_1,d_0$ defined in \cite{WW17} up to a constant.
As a consequence, a series of quarternionic pluripotential results associated
to theses two operators \cite{Wan20a}\cite{Wan20b}\cite{Wan17}\cite{Wan19}\cite{WW17}
\cite{WK17}\cite{WZ15} can be used without hesitation in terms of the intrinsic operators $\p,\pj$.
Based on Bedford-Taylor's theory, Wan and Wang \cite{WW17} proved that quaternionic Monge-Amp\`ere operator can
be defined for locally bounded quaternionic plurisubharmonic functions and 
is continuous in decreasing sequence. Motived by Cegrell's work \cite{Ceg98}
\cite{Ceg04}\cite{Lu15}, the definition can be extended to some classes of unbounded functions \cite{Wan20a}
More precisely, Wan \cite{Wan20a} showed that the quaternionic Monge-Amp\`ere
operator is well defined in the Cegrell type classes $\mathcal{E}(\Omega)$ and 
$\mathcal{E}_p(\Omega),p\geq1$. Some energy estimates were also established and by using variational method in \cite{Lu15}
she proved that a Radon measure $\mu=(\pp \varphi)^n$ with 
$\varphi\in\mathcal{E}_p(\Omega)$ if and only if $\mathcal{E}_p(\Omega)\subset L^p(\Omega),p\geq1.$
When $0<p<1$, the variational method becomes invalid since the constant of p-energy estimate, which plays an important
role in the variational approach, exceeds 1 in an example shown by Do and Nguyen \cite{DN25}. 

However, we can
overcome this difficulty by using ideas from \cite{Sal23}\cite{Ben15} and solve the quaternionic
Monge-Amp\`ere equation for a more general Cegrell type class $\mathcal{E}_{\chi}(\Omega)$
where $\chi$ is a weight function (see the definition in Section \ref{weighted energy classes}).
When $\chi=-(-t)^p,p>0$ and
$\mathcal{E}_{\chi}(\Omega)$ coincides with $\mathcal{E}_{p}(\Omega)$, while $p$ is not
necessarily bigger than 1. Such class is initiated in \cite{GZ07} and has led to
extensive subsequent researches in the decade (see e.g. \cite{DDL23} and the references therein).  As a first step to study the equation, we prove that $\mathcal{E}_p(\Omega)\subset\mathcal{E}(\Omega),p>0$(Theorem \ref{Ep in E})
and furthermore,
\begin{theorem}(Proposition \ref{energy class in cegrell class})
    $$\mathcal{E}_{\chi}(\Omega)\subset\mathcal{E}(\Omega).$$
    In particular, the quaternionic Monge-Amp\`ere operator $(\pp u)^n$ is well-defined
    for any $u\in\mathcal{E}_{\chi}(\Omega)$ and $-\chi(u)\in L^1((\pp u)^n)$. Moreover, $(\pp u)^n$ puts no mass on quaternionic polar subsets.
\end{theorem}

The proof is inspired by \cite{Ben11}\cite{BGZ08}.
In the complex setting, the energy
classes $\mathcal{E}_{\chi}(\Omega)$ is studied in many papers, to cite a few, e.g. \cite{Ceg98}\cite{Ceg04}\cite{Ben09}\cite{Ben13}
\cite{DDL23}
\cite{DDNL18}\cite{Sal23}.

An interesting question in analogy with a conjecture of Guedj-Zeriahi in \cite{GZ07} is how to characterize finite weighted energy range of 
quaternionic Monge-Amp\`ere operator. In the complex setting, a general solution of Guedj-Zeriahi's 
conjecture in compact K\"ahler manifolds is obtained in \cite{DDL23}, while in the local context
the characterization is also available, see \cite{Sal23}\cite{Ben11}\cite{Ben13}\cite{Ben15} for different
settings on the weight function $\chi$. Our second main theorem is to give a characterization of
finite enery range of \emph{quaternionic} Monge-Amp\`ere operator:
\begin{theorem}(Theorem \ref{characterization})
    Let $\mu$ be a positive Radon measure, and let $\chi \in \mathcal{W}^{-}\cup\mathcal{W}_M^{+}$. The following conditions are equivalent:
    \begin{enumerate}
        \item there exists a \textit{unique} function $\varphi \in \mathcal{E}_{\chi}(\Omega)$ such that $\mu = (dd^c \varphi)^n$;
        \item\label{cond2} $\chi(\mathcal{E}_{\chi}(\Omega)) \subset L^1(d\mu)$;
        \item\label{cond3} there exists a constant $C > 0$ such that
        \[
        \int_{\Omega} -\chi \circ \psi \, d\mu \leq C,
        \]
        for all $\psi \in \mathcal{E}_0(\Omega)$, $E_{\chi}(\psi) \leq 1$;
        \item\label{cond4} there exists a positive constant $A$ such that
        \[
        \int_{\Omega} -\chi \circ \psi \, d\mu \leq A \max(1, E_{\chi}(\psi)), \quad \forall \psi \in \mathcal{E}_0(\Omega).
        \]
    \end{enumerate}

\end{theorem}
Here we let $\mathcal{W}^{-}$ (resp $\mathcal{W}^{+}$) be the set of convex (resp concave) increasing functions
$\chi:\mathbb{R}^{-}\to\mathbb{R}^{-}$ such that $\chi(-\infty)=-\infty$ and $\chi(0)=0$ and  $\mathcal{W}_M^{+}$ be the set of  
functions $\chi\in\mathcal{W}_M^{+}$ with the property
$$|t\chi\prime (t)|\leq M|\chi(t)|,\forall t\in\mathbb{R}^{-}.$$

As an important step
to the proof, we show that fine property (Proposition \ref{plurifine property}) holds in $\mathcal{E}(\Omega)$.
As a corollary, we prove a mass concentration theorem for the quaternionic plurisubharmonic
envelope $P(f):=(\sup\{\varphi\in \mpsh:\varphi\leq f\})^{*}$ where $f$ is a measurable function on $\Omega$.  
\begin{theorem}(Theorem \ref{mass conc})
    Assume that $f$ is quasi-continuous and $f\leq0$, and there is 
        $\psi\in\mathcal{E}(\Omega)$ such that $\psi\leq f$. Then the measure $(\p\pj P(f))^n$ is concentrated
    on $\{P(f)=f\}$.
\end{theorem}
When $u$ is continous and $\psi\in\mathcal{E}_1$, the theorem is proved in \cite{Wan20a}.
For more applications of \emph{complex} plurisubharmonic envelope, the reader can refer
to \cite{DDL23}\cite{BBGZ13}\cite{LN22}
\cite{Lu15}\cite{GLZ19}\cite{Sal23}.

We should remark that the assumption of quaternionic strictly psedoconvexity on $\Omega$
can be relaxed to quaternionic hyperconvexity as in \cite{Wan20a} if we 
do not consider Cegrell type class with sero smallest maximal majorant, namely $\mathcal{N}(\Omega)$
whose definition requires quaternionic strictly psedoconvexity.

\section*{acknowledgement}
 I wish to express my gratitude to my mentors Stephan S.-Yau and Shing-Tung Yau for 
 their guidance and support. I also wish to  thank Mohammed Salouf for some helpful discussions and explanations on
some details in his paper \cite{Sal23} sincerely. 
\section{Preliminary}
For general references for quarternionic linear algebra and basic properties of quaternionic
plurisubharmonic functions, one can refer to \cite{Ale03a}\cite{Ale03b}\cite{AV10}
which initiated the study of quarternionic Monge-Amp\`ere equation.

Based on Bedford-Taylor's theory for degenerate complex Monge-Amp\`ere equation, the authors
in \cite{Wan17}\cite{Wan19}\cite{WW17}\cite{WK17}\cite{WZ15} reformulated  classical theorems of
pluripotential theory into quarternionic ones. Thesis of Sroka \cite{Sro21} is also a 
good exposition for this topic. Here we only give a brief introduction.

Let us fix a quarternionic algebra $$\mathbb{H}^n=\{x_0+x_1\mathbf{i}+x_2\mathbf{j}+x_3\mathbf{k}|x_0,x_1,x_3\in\mathbb{R}\},$$
where $\mathbf{i},\mathbf{j},\mathbf{k}$ satisfy quaternionic relations and we can consider $\mathbb{H}^n$ as a natural
right module. The right multiplications $i,j,k$ induce three complex structures on $\mathbb{H}^n$ 
when treating $\mathbb{H}^n$ as a flat hypercomplex manifold. Let us denote them by 
$I,J,K$. $\mathbb{H}^n$ has two natural coordinate systems, among which $\mathbb{H}^n$
can be considered as $\mathbb{C}^{2n}$ by
$$(q_0,q_1,...,q_{n-1})\in \mathbb{H}^n\mapsto (z_0,z_1,...,z_{2n-1})\in\mathbb{C}^n$$
where $q_i=z_{2i}+\mathbf{j}z_{2i+1},i=0,...,n-1.$ Another coordinate is  given by 
$$(q_0,q_1,...,q_{n-1})\in\mathbb{H}^n\mapsto (x_0,x_1,...,x_{4n-1})\in\mathbb{R}^{4n}$$
where $q_i=x_{4i}+x_{4i+1}\mathbf{i}+x_{4i+2}\mathbf{j}+x_{4i+3}\mathbf{k}$.
It is easy to see that $z_j=x_{2j}+(-1)^jx_{2j+1}\mathbf{i}$ for $j=0,...,2n-1$. 

Let $\partial,\bar\partial$ be the canonical differential operators induced by the complex
structure $I$ and $d=\partial+\bar\partial,d^c=\mathbf{i}(\bar\p-\p)$. The twisted differential
is introduced in \cite{AV10}\cite{Ver02} as
$$\pp:=J^{-1}\circ\bar\p\circ J.$$
This plays the role of $\bar\p$ in the hypercomplex setting. The interested readers 
can refer to the papers mentioned above. 
From a 
basic computation, we have explicit expressions for the hypercomplex structures in terms of
natural frames of $\mathbb{H}^n$:
\begin{proposition}\cite{Sro21}
    The induced complex structures $I,J,K$ are given in the frame $\{\p_{x_i}\}$ or $\{\p_{z_j}\}$ by
    \begin{equation*}
        \begin{aligned}
            &(\p_{x_{4i}})I=\p_{x_{4i+1}},(\p_{x_{4i+2}})I=-\p_{x_{4i+3}} \\
            &(\p_{x_{4i}})J=\p_{x_{4i+2}},(\p_{x_{4i+1}})J=\p_{x_{4i+3}}\\
            &(\p_{x_{4i}})K=\p_{x_{4i+1}},(\p_{x_{4i+1}})K=-\p_{x_{4i+2}},\\
            &(\p_{z_{2i+1}})J=\p_{\overline{z_{2i+1}}},(\p_{z_{2i+1}})J=-\p_{\overline{z_{2i}}}.
        \end{aligned}  
    \end{equation*}
    For a smooth function \( u : \mathbb{H}^n \to \mathbb{R} \) we have 
\begin{equation*}
    \begin{aligned}
\partial_J u &= (J^{-1} \bar{\partial} J) u = J^{-1} (\bar{\partial} u) = J^{-1} \left( \sum_{j=0}^{2n-1} \p\overline{{z_j}} u d\bar{z}_j \right)\\
&= J^{-1} \left( \sum_{j=0}^{2n-1} \p_{\overline{{z_{j+(-1)^j}}}} u d{\overline{z_{{j+(-1)^j}}}} \right) ,\\
 &= \sum_{j=0}^{2n-1} \p_{\overline{{z_{j+(-1)^j}}}} u (-1)^{j+1} dz_j\\
 \partial \partial_J u &= \sum_{i,j} \left( (-1)^{j+1} \partial_{z_i} \p_{\overline{{z_{j+(-1)^j}}}} u \right) dz_i \wedge dz_j\\
&= \sum_{i<j} \left( (-1)^{j+1} \partial_{z_i} \p_{\overline{{z_{j+(-1)^j}}}} u - (-1)^{i+1} \partial_{z_j} \p_{\overline{{z_{i+(-1)^i}}}} u \right) dz_i \wedge dz_j.
\end{aligned}
\end{equation*}
Especially
$$\partial \partial_J \left( \frac{1}{2} \sum_{k=0}^{2n-1} z_k \overline{z_k} \right) = \sum_{k=0}^{n-1} dz_{2k} \wedge dz_{2k+1}:=\Omega.$$
\end{proposition}
Let $\Lambda_I^{p,q}(\mathbb{H}^n)$ be the space of $(p,q)$-forms with respect to $I$, then 
$$\pj:\Lambda_I^{k,0}(\mathbb{H}^n)\to\Lambda_I^{k+1,0}(\mathbb{H}^n)$$
satisfies $$\p\pj+\pj\p=0,$$
$$\pj^2=0.$$
Suppose that $f:\mathbb{H}^n\to\mathbb{H}$ is a $C^2$ function. The formal quaternionic
derivative is defined as 
\begin{equation*}
    \begin{aligned}
        \frac{\p f}{\p\bar{q}_{\alpha}}&=\frac{\p f}{\p x_{4\alpha}}+\mathbf{i}\frac{\p f}{\p x_{4\alpha+1}}+\mathbf{j}\frac{\p f}{\p x_{4\alpha+2}}+\mathbf{k}\frac{\p f}{\p x_{4\alpha+3}},\\
        \frac{\p f}{\p q_{\alpha}}&=\frac{\overline{p \bar{f}}}{\p \bar{q}_{\alpha}}=\frac{\p f}{\p x_{4\alpha}}-\frac{\p f}{\p x_{4\alpha+1}}\mathbf{i}-\frac{\p f}{\p x_{4\alpha+2}}\mathbf{j}-\frac{\p f}{\p x_{4\alpha+3}}\mathbf{k}.
    \end{aligned}
\end{equation*}
It is easy to verify that the matrix
$$Hess(f,\mathbb{H})=\left(\frac{\p^2 f}{\p\bar{q_{\alpha}}\p q_{\beta}}\right)_{\alpha,\beta\in\{1,...,n\}}$$
is hyperhermitian for any real-valued $f$. 
Using the Moore determinant of a hyperhermitian matrix we can define the quaternionic Monge-Amp\`ere operator for a smooth
function $u$:
$$(\pp u)^n=\frac{n!}{4^n}det\left(\frac{\p^2 u}{\p_{{\bar{q}}_l}\p_{q_k}}\right)dz_0\wedge dz_1\wedge...\wedge dz_{2n-2}\wedge dz_{2n-1},$$
see \cite{AV06}. See also \cite{WW17} in terms of a different setting. It is observed by Sroka that 
the operators $d_0,d_1$ in \cite{WW17} have the following relations with $\p,\pj$:
\begin{proposition}\cite{Sro20}
    For the basis $\omega^k=(-1)^kdz_{k+(-1)^k}$,
    $$d_0=2\pj,d_1=-2\p, ~and~\Delta=d_0d_1=4\p\pj.$$
\end{proposition}
From above we are able to use all results from \cite{Wan20a}\cite{Wan17}\cite{Wan19}\cite{WW17}\cite{WK17}
\cite{WZ15}. Recall the quarternionic plurisubharmonic functions defined in \cite{Ale03b}:
\begin{definition}
    Let $\Omega\in\mathbb{H}^n$ be a domain. An upper semicontinuous function $f:\Omega\to
    \mathbb{R}$ quaternionic plurisubharmonic, qpsh for short, if for any affine 
    quaternionic line intersected with $\Omega$, $f$ is subharmonic  when restricted
    in such line. The set of all qpsh functions on $\Omega$ is denoted by $QPSH(\Omega)$.
\end{definition}
It can be shown that a smooth function $u$ being qpsh if and only if $\pp u\geq0$. 
The notions of positivity on hyperhermitian manifolds are introduced in \cite{AV06}. The reader can
also refer to \cite{Sro21}. On the other hand, classic results of complex pluripotential
theory (see \cite{BT82}) have  analogues in the quarternionic setting due to e.g.
\cite{WW17}\cite{WZ15} where the authors studied locally bounded qpsh functions. A 
comparison principle for locally bounded qpsh function is proved in \cite{WZ15}.
\begin{theorem}\cite{WZ15}
    Let $u,v\in\mpsh\cap L_{loc}^{\infty}(\Omega).$ If $\liminf_{q\to\p\Omega}(u(q)-v(q))\geq0$,
    then 
    $$\int_{\{u<v\}}(\pp v)^n\leq\int_{\{u<v\}}(\pp u)^n.$$
    In particular, if $(\pp v)^n\geq (\pp u)^n$ as measures, then $u\geq v$ in $\Omega$.
\end{theorem}

\section{Cegrell type classes}\label{S1}
In \cite{Wan20a} the author introduced and studied the following 
Cegrell type classes of QPSH functions in $\Omega$. Let $p\geq0$. In the rest of 
this paper, we fix $\Omega$ to be a quaternionic strictly pseudoconvex domain defined as below.
\begin{definition}
    A smoothly bounded domain $\Omega\subset\mathbb{H}^n$ is quaternionic strictly pseudoconvex if
    there exists a domain $U$ and a function $\rho$ satisfying
    \begin{equation*}
        \begin{cases}
            &\rho\in\mpsh\cap C^2(\Omega)\\
            & D\subset\subset U\\
            & D=\{\rho<0\}\\
            &d\rho\neq 0~on~\p\Omega\\
            &(\pp \rho)^n\geq \Omega_n ~on~U
        \end{cases}
    \end{equation*}
    where $\Omega_n:=\frac{\Omega^n}{n!}=dz_0\wedge dz_1\wedge...\wedge dz_{2n-2}\wedge dz_{2n-1}$.
\end{definition}
\begin{definition}
    \begin{equation*}    
        \begin{split}
    &\mathcal{E}_0(\Omega)=\{u\in\mpshn\cap L^{\infty}(\Omega):u=0~on~\p\Omega
    ~and~\int_{\Omega}(\p\pj u)^n<+\infty\};\\
    & \mathcal{E}(\Omega)=\{u\in\mpshn:\forall z\in\Omega,\exists
    V~a~neighborhood~of~z,\exists(u_j)_j\subset\mathcal{E}_0(\Omega),
    u_j\downarrow u ~on~V\\
    &~and~\sup_j\int_{\Omega}(\p\pj u_j)^n<+\infty\};\\
    &\mathcal{F}(\Omega)=\{u\in\mpshn:\exists u_j\in\mathcal{E}_0(\Omega),
    u_j\downarrow u~and~\sup_j\int_{\Omega}(\p\pj u_j)^n<+\infty\};\\
    &\mathcal{E}_p(\Omega)=\{u\in\mpsh:\exists(u_j)_j\in\mathcal{E}_0(\Omega),u_j\downarrow u~and~\sup_j\int_{\Omega}|u_j|^p(\p\pj u_j)^n<+\infty\}.
        \end{split}
    \end{equation*}
\end{definition}
\begin{theorem}\label{Ep in E}
As in the complex case, these classes satisfy the inclusion
$$ \mathcal{E}_0(\Omega)\subset\mathcal{E}_p(\Omega)\cap\mathcal{F}(\Omega)\subset\mathcal{E}_p(\Omega)\cup
\mathcal{F}(\Omega)\subset\mathcal{E}(\Omega).$$
\end{theorem}
\begin{proof}
    The first and second inclusion follow from the definitions. For the third, it suffices to 
    prove that $\mathcal{E}_p(\Omega)\subset\mathcal{E}(\Omega)$.

    Fix $u \in\mathcal{E}_{p}(\Omega)$.
    Assume that $\{u_j\}\in\mathcal{E}_0$ decreases to $u$ on $\Omega$ and satisfies  $\sup_j\int_{\Omega}(-u_j)^p(\pp u_j)^n<+\infty$.
    Take $\omega\subset\subset\Omega$ and set $u_{j\omega}:=\sup\{v\in\mpshn:v\leq u_j~on~\omega\}$.
    Then we have $u_j\leq u_{j\omega},u_{j\omega}\in\Omega_0(\Omega)$ and $\supp (\p\pj u_{j\omega})\in\overline{\omega}$.

    If $0<p<1$, integration by parts gives
    $$\int_{\Omega}(-u_j)^p(\pp u_{j\omega})^n\leq\sup\int_{\Omega}(-u_j)^p(\pp u_j)^n.$$
    Taking the supreme, we have
   \begin{equation}
        \begin{aligned}
                \int_{\Omega}(\pp u_{j\omega})^n&\leq\frac{1}{\inf_{\omega}(-u_1)}\int_{\Omega}(-u_1)^p(\pp u_{j\omega})^n\\
            &\leq \sup_j\frac{1}{\inf_{\omega}(-u_1)}\int_{\Omega}(-u_j)^p(\pp u_{j\omega})^n\leq\infty.
        \end{aligned}
    \end{equation}
    If $p\geq1$, since $u_j\leq u_{j\omega}$, by \cite[Corollary 3.4]{Wan20a} there exists a constant $C>0$
    such that
$$\int_{\Omega}(-u_{j\omega})^p(\pp u_{j\omega})^n\leq C\int_{\Omega}(-u_j)^p(\pp u_j)^n.$$
Taking the supreme, we have
\begin{equation*}
    \begin{aligned}
        \int_{\Omega}(\pp u_{j\omega})^n&\leq\frac{1}{\inf_{\omega}(-u_{1\omega})}\int_{\Omega}(-u_{j\omega})^p(\pp u_{j\omega})^n\\
        &<\infty.
    \end{aligned}
\end{equation*}
Hence for $p>0$, $\mathcal{E}_p(\Omega)\subset\mathcal{E}(\Omega)$.
\end{proof}

If $u\in\mathcal{E}(\Omega)$, then $(\p\pj u)^n$ defines a positive
Radon measure by \cite[Theorem 3.2]{Wan20a}. The set $\mathcal{E}(\Omega)$
is the largest set for which the quaternionic Monge-Amp\`ere operator
$(\p\pj\cdot)^n$ is well-defined and continuous along decreasing sequence in a 
certain sense in \cite[Corollary 3.3]{Wan20a}, which is proved by a similar argument of \cite{Ceg04}.

Let $\Omega_j$ be a fundamental sequence of quaternionic
strictly pseudoconvex subdomains of $\Omega$ and 
$u\in\mathcal{E}(\Omega)$, we define 
\begin{equation*}
    u^j=\sup\{\varphi\in\mpshn:\varphi\leq u~on~\Omega_j^c\}.
\end{equation*}
Note that we have $u\leq u^j\leq u^{j+1}$ by definition and hence
$u^j\in\mathcal{E}(\Omega)$ for any $j$, and so is $(\tilde{u}=\lim_{j\to\infty}u^j)^{*}$
 the smallest upper semicontinuous majorant of $\lim_j u^j$. Recall that a function $u\in\mpsh$
 is maximal if for any $v\in\mpsh$ and compactly open $U\subset\Omega$, $u\geq v$ holds on $U$ provided
 $u\geq v$ on $\p U$.
 \begin{proposition}
    The definition of $\tilde{u}$ is independent of the sequence $\Omega_j$ and $\tilde{u}\in\mathcal{E}(\Omega)$
    is maximal.
 \end{proposition}
\begin{proof}
    We only need to prove the second statement. Approximate $u$ by a decreasing sequence $u_k\in\mathcal{E}_0\cap C^0(\Omega)$and
    $u_k=0$ in $\p\Omega$ by \cite[Theorem 3.1]{Wan20a}.
   Let $\Psi_k^j$ be the unique solution to the following Dirichlet problem
   \begin{equation*}
    \begin{cases}
        &v\in\mpsh\cap C(\overline{\Omega_j})\\
        &(\p\pj v)^n=0\\
        &v|_{\p\Omega}=u_k\in C(\p\Omega_j).
    \end{cases}
   \end{equation*}
   By a theorem of Alesker \cite{Ale03a} we have that such solution exists. It follows 
   from \cite[Corollary 5.1]{WZ15} that $\Psi_k^j$ is maximal. We can extend the definition of $\Psi_k^j$
   on the whole $\Omega$ by define it to be $u_k$ in $\Omega_j^c$. Still denote it by $\Psi_k^j$ and we have $\Psi_k^j\in\mpshn$.
    Since $u_k$ is decreasing to $u$,
   it follows from the quaternionic analogue of Perron-Bremermann construction that
   $\Psi_k^j$ is decreasing with respect to $k$. We claim that the limit $\Psi_j:=\lim_k\Psi_k^j$ is
   equal to $u^j$.

   On the one hand, it is easy to see that $\Psi^j\in\mpshn\leq u$ on $\Omega_j$ and hence it lies in the defining family
   of $u^j$. So we have $\Psi^j\leq u^j$. On the other hand, $u^j\leq\Psi_k^{j+1}$ on $\p\Omega_j$ since 
   $u^j$ is equal to $u$ on $\p\Omega_j$ by definition. Since $\Omega_j\subset\subset\Omega_{j+1}$
   and $\Psi_k^{j+1}$ is maximal on $\Omega_{j+1}$, it follows from the maximality that $u^j\leq\Psi_k^{j+1}$ on $\Omega_j$.
   Hence by taking the limit $u^j\leq\Psi^j$. The claim holds.
   
   Applying the monotonically convergence theorem of quaternionic Monge-Amp\`ere operator, we have 
   that $u^j$ is maximal in $\Omega_j$ and $\tilde{u}$ is hence maximal on $\Omega$.

\end{proof}
If $u\in\mathcal{F}(\Omega)$ or $\mathcal{E}_p(\Omega)$, then so is  $\tilde{u}$ \cite[Corollary 3.3]{Wan20a} and equals to zero (e.g.see 
\cite[Lemma 4.7]{Wan20a} and \cite[Theorem 1.3]{WZ15}).
Denote the class $\mathcal{N}(\Omega)$ to be the set of $u\in\mathcal{E}(\Omega)$ such 
that $\tilde{u}:=(\lim u_j)^{*}=0$. We have $\mathcal{F}(\Omega)\subset\mathcal{N}(\Omega)$
and $\mathcal{E}_p(\Omega)\subset\mathcal{N}(\Omega),\forall p$. We also denote by $\mathcal{E}^a(\Omega)$
the set of $u\in\mathcal{E}(\Omega)$ such that $(\pp u)^n$ puts no mass on quaternionic polar subsets.
\begin{comment}
\begin{theorem}\label{FEcp}
    If $u\in\mathcal{F}(\Omega)$ and $v\in\mathcal{E}(\Omega)$, then
    \begin{equation}
        \int_{\{u<v\}}(\p\pj v)^n\leq \int_{\{u<v\}\cup\{u=-\infty\}}(\p\pj u)^n.
    \end{equation}
   
\end{theorem}
\begin{proof}
    The proof follows from a analogous argument in \cite[Corollary 3.6]{Ceg08}.
\end{proof}
\begin{theorem}
    Suppose $u\in\mathcal{F}^a(\Omega),v\in\mathcal{E}(\Omega)$ and $(\p\pj u)^n\leq(\p\pj v)^n$.
    Then $v\leq u$ on $\Omega$.
\end{theorem}
\begin{proof}
    Via Theorem\ref{FEcp}, the proof of \cite[Corollary 1.1]{WZ15} can be refined to this case.
\end{proof}
\end{comment}

Now we will mainly prove a mass concentration theorem for 
quaternionic plurisubharmonic envelop as in the complex case.
In the complex case, plurisubharmonic envelops plays an important role in
the study of complex Monge-Amp\`ere equation (see \cite{Sal23} and a general paper \cite{DDL23} and 
references therein). Here on the quarternionic domain  we can 
we can define a  similar envelop associated to a measurable function
$f$ by
$$P(f):=(\sup\{\varphi\in \mpsh:\varphi\leq f\})^{*}.$$

First we have the following proposition.

\begin{proposition}
If $f$ is a measurable function on $\Omega$, then 
$$P(f)=\sup\{\varphi\in\mpsh:\varphi\leq f ~\mbox{quasi-everywhere}\},$$
where the term quasi-everywhere means outside a pluripolar set. Moreover,
if $(f_j)$ is a decreasing sequence of measurable functions 
converging to $f$, then $P(f_j)$ decreases to $P(f)$.
\end{proposition}
\begin{proof}
    By relying on the quaternionic potential theory originated from Bedford-Taylor
    and reproved by Wan-Kang \cite{WK17}, one can obtain a proof 
    analogous to \cite{Sal23}.
    
\end{proof}
Recall the definition of the quarternionic Monge-Amp\`ere capacity: Given a Borel subset $E\subset\Omega$,
we set $$ Cap(E):=\sup\bigg\{\int_E(\p\pj u)^n:u\in\mpsh, -1\leq u\leq0
    \bigg\}.$$ Then we have 
    \begin{lemma}\label{cap estimate}
        Let $u\in\mathcal{E}_0(\Omega)$. Then for all $s>0$ and $t>0$, we have
        $$t^n Cap(\{u<-s-t\})\leq\int_{\{u<-s\}}(\pp u)^n\leq s^n Cap(\{u<-s\}).$$
    \end{lemma}

    A function $f$ is quasi-continuous if for every $\epsilon>0$, 
    there is a Borel set $E\subset\Omega$ such that $Cap(E)\leq\epsilon$
    and $f$ is continuous on $\Omega\setminus E$. 

    Now we can prove the mass concentration theorem about the QPSH envelope:
    \begin{theorem}\label{mass conc}
        Assume that $f$ is quasi-continuous and $f\leq0$, and there is 
        $\psi\in\mathcal{E}^{a}(\Omega)$ such that $\psi\leq f$. Then 
        $P(f)\in\mathcal{E}^{a}(\Omega)$ and $(\p\pj P(f))^n$ is concentrated
        on $\{P(f)=f\}$.
    \end{theorem}
    
\begin{proof}
    First we note that the hypothesis $\psi\leq f$ for some $\psi\in\mathcal{E}^a$ guarantees that $P(f)\in\mathcal{E}^a.$
    Indeed, let $P$ be a compact Q-polar subset and $h\in\mathcal{E}_0(\Omega)$ such that $h\leq -1$
    near $P$. Then we have a sequence $0\geq h_p\geq h$ such that $\{h_p\}_p\in\mathcal{E}_0$ equals to $-1$ on $P$ and
    increases to $0$ outside a Q-polar subset \cite[Proposition 3.4]{WK17}.  Integration by parts
    gives that $\int_{\Omega}-h_p(\pp P(f))^n\leq\int_{\Omega}-h_p(\pp \psi)^n$. Since $(\pp\psi)^n$
    is non-Q-polar, the first statement holds by letting $p$ tend to $+\infty$.

    For the second statement,
    the proof is a quarternionic analogue to that of \cite[Theorem 2.4]{Sal23}. To treat the case that $f$ is not
    bounded from below, we can approximate it by $f_j:=\max(f,-j)$ and use the result in which
    $f$ is bounded from below. As an important step to treat this general case, one need the following
    fine property (Proposition \ref{plurifine property} )for functions in $\mathcal{E}(\Omega)$.

\end{proof}
\begin{proposition}\label{plurifine property}
    Let $u,u_1,...,u_{n-1}\in\mathcal{E}(\Omega),v\in\mpshn$ and $T=\pp u_1\wedge...\wedge\pp u_{n-1}$.
    Then $$\pp\max(u,v)\wedge T|_{u>v}=\pp u\wedge T|_{u>v}.$$
    In particular, the sequence $\mu^j:=\mathbf{1}_{\{u>-j\}}(\pp\max(u,-j))^n$ is an increasing sequence of
    positive measures. We denote by $\mu_{u}$ its limit.
\end{proposition}
\begin{proof}
     The proof is divided into two steps.

     a) First we prove the proposition for $v \equiv a < 0$. Assume that $u, u_1, ..., u_{n-1} \in \mathcal{F}$. Using Theorem 3.1 in \cite{Wan20a} we can find
    \[
    \mathcal{E}_0 \cap C(\bar{\Omega}) \ni u^j \searrow u, \quad \mathcal{E}_0 \cap C(\bar{\Omega}) \ni u_k^j \searrow u_k, \quad k = 1, ..., n - 1.
    \]
    Since $\{u^j > a\}$ is open we have
    \[
    \pp \max(u^j, a) \wedge T_j|_{\{u^j>a\}} = \pp u^j \wedge T_j|_{\{u^j>a\}}.
    \]
    Thus from the inclusion $\{u > a\} \subset \{u^j > a\}$ we obtain
    \[
    \pp \max(u^j, a) \wedge T_j|_{\{u>a\}} = \pp u^j \wedge T_j|_{\{u>a\}},
    \]
    where $T_j = \pp u_1^j \wedge ... \wedge \pp u_{n-1}^j$. It follows from the following Proposition \ref{11} that
    \[
    \max(u - a, 0) \pp \max(u^j, a) \wedge T_j \to \max(u - a, 0) \pp \max(u, a) \wedge T,
    \]
    \[
    \max(u - a, 0) \pp u^j \wedge T_j \to \max(u - a, 0) \pp u \wedge T.
    \]
    Hence
    \[
    \max(u - a, 0)[\pp \max(u, a) \wedge T - \pp u \wedge T] = 0.
    \]
    Using a basis fact (see \cite[Lemma 4.2]{KH09}) we have
    \[
    \pp \max(u, a) \wedge T = \pp u \wedge T \quad \text{on } \{u > a\}.
    \]
    b) Assume that $v \in QPSH^-(\Omega)$. Since $\{u > v\} = \bigcup_{a \in \mathbb{Q}^-} \{u > a > v\}$, it suffices to show that
    \[
    \pp \max(u, v) \wedge T = \pp u \wedge T \quad \text{on } \{u > a > v\}
    \]
    for all $a \in \mathbb{Q}^-$, the set of negative rational numbers. Since $\max(u,v) \in \mathcal{E}$, by a) we have
    \begin{align}
    \pp \max(u,v) \wedge T|_{\{\max(u,v)>a\}} &= \pp \max(\max(u,v), a) \wedge T|_{\{\max(u,v)>a\}}, \\
    &= \pp \max(u, a) \wedge T|_{\{\max(u,v)>a\}}, \\
    \pp u \wedge T|_{\{u>a\}} &= \pp \max(u, a) \wedge T|_{\{u>a\}}.
    \end{align}
    Since $\max(u,v,a) = \max(u,a)$ on the open set $\{a > v\}$, we have
\begin{equation}
    \pp \max(u,v,a) \wedge T|_{\{a>v\}} = \pp \max(u,a) \wedge T|_{\{a>v\}}.
\end{equation}

Since $\{u > a > v\} \subset \{u > a\}, \{a > v\}, \{\max(u,v) > a\}$ and (2.1), (2.2), (2.3) we have
\[
    \pp \max(u,v) \wedge T|_{\{u>a>v\}} = \pp u \wedge T|_{\{u>a>v\}}.
\]

\end{proof}
The following two results are derived from the above proposition.
\begin{corollary}\label{coro 5.4}
    Let $u,v \in \mathcal{E}(\Omega)$, and let $\mu$ be a positive Radon measure vanishing on quaternionic polar sets. If
\begin{equation*}
    ( \pp u)^n \geq \mu \quad \text{and} \quad ( \pp v)^n \geq \mu,
\end{equation*}
then
\begin{equation*}
    ( \pp \max(u,v))^n \geq \mu.
\end{equation*}
\end{corollary}
\begin{proof}
    By Proposition \ref{plurifine property}, we have
\begin{equation*}
    ( \pp \max(u,v))^n \geq \mathbf{1}_{\{u > v\}}(  \pp u)^n + \mathbf{1}_{\{u < v\}}( \pp v)^n \geq \mathbf{1}_{\{u \neq v\}} \mu.
\end{equation*}
If $\mu(\{u = v\}) = 0$, then we are done. From a similar argument as \cite[Corollary 1.10]{GZ07} we can derive that 
\begin{equation*}
    \mu(\{u = v + t\}) = 0, \quad \forall t \in \mathbb{R} \setminus I,
\end{equation*}
where $I$ is at most countable. Take $\varepsilon_j \in \mathbb{R} \setminus I$, $\varepsilon_j \searrow 0$. We have
\begin{equation*}
    (\pp \max(u, v + \varepsilon_j))^n \geq \mu.
\end{equation*}
Let $\varepsilon_j \to 0$ and we can complete the proof. 
\end{proof}
\begin{corollary}\label{coro 5.5}
    Let $u,v \in \mathcal{E}(\Omega)$. If $( \pp v)^n$ vanishes on pluripolar sets, then
\begin{equation*}
    ( \pp \max(u,v))^n \geq 1_{\{u > v\}}( \pp u)^n + 1_{\{u \leq v\}}(\pp v)^n.
\end{equation*}
In particular, if moreover $u \geq v$, then
\begin{equation*}
    \mathbf{1}_{\{u = v\}}(\pp u)^n \geq \mathbf{1}_{\{u = v\}}( \pp v)^n.
\end{equation*}
\end{corollary}
\begin{proof}
    The first statement follows from Corollary \ref{coro 5.4} and the second 
    follows from the first.
\end{proof}

\begin{proposition}\label{11}
    Suppose that $u^k\in\mathcal{F}(\Omega),k=1,...,n$ and $\{g_k^j\}\subset\mathcal{E}_0(\Omega)$
    decreases to $u^k,j\to +\infty.$
    Then for each $h\in\mpshn$ we have
    $$\lim_{j\to +\infty}\int_{\Omega}h\pp g_j^1\wedge...\wedge g_j^n=\int_{\Omega}h\pp u^1\wedge...
    \wedge\pp u^n.$$
\end{proposition}
\begin{proof}
    
    The proof is obtained in \cite[Theorem 3.2]{Wan20a} under the assumption that
    $\sup_j\int_{\Omega}\pp g_j^1\wedge...\wedge \pp g_j^n<+\infty$. However, this assumption is 
    vacuous. Indeed, since $\mathcal{F}(\Omega)\subset\mathcal{E}(\Omega)$, by \cite[Corollary 3.3]{Wan20a}
    we know that $\pp g_j^1\wedge...\wedge \pp g_j^n$ weakly convergent to $u^1,...,u^n$. Since $\Omega$ 
    is open, $\limsup_j \int_{\Omega}\pp g_j^1\wedge...\wedge \pp g_j^n\leq\int_{\Omega}\pp u^1\wedge...\wedge\pp u^n$ which is finite.

\end{proof}
\begin{remark}
    The above proposition is a quarternionic analoge of \cite[Proposition 5.1]{Ceg04}.
\end{remark}
\section{Weighted energy classes}\label{weighted energy classes}
In this section we study the existence of solutions to the following Dirichlet problem
\begin{equation}
    (\pp u)^n=\mu, u\in\mathcal{E}_{\chi},
\end{equation}
where $\mu$ is a positive Radon measure vanishing on Q-polar subsets.
The motivation comes from \cite{Sal23}\cite{Ben09}\cite{Ben15} and it is desired to give a characterization of finite energy range of
quarternionic Monge-Amp\`ere operator on quaternionic domain $\Omega\subset \mathbb{H}^n$. This will generalize the result in which the special
concave weight $\chi(t)=-(-t)^p(p\geq1)$ is considered.

Recall that  $\mathcal{W}^{-}$ (resp$\mathcal{W}^{+}$) is the set of convex (resp concave) increasing functions
$\chi:\mathbb{R}^{-}\to\mathbb{R}^{-}$ such that $\chi(-\infty)=-\infty$ and $\chi(0)=0$ and the set $\mathcal{W}_M^{+}$ consists of 
functions $\chi\in\mathcal{W}_M^{+}$ with the property
$$|t\chi\prime (t)|\leq M|\chi(t)|,\forall t\in\mathbb{R}^{-}.$$
By \cite{TV21} in this case we have 
\begin{equation}\label{subh ineq}
    -\chi(ct)\leq -c^M\chi(t)~ for~ all~ t\leq0~ and~ all ~c\geq1.
\end{equation}
 Thus we call 
$\chi\in\mathcal{W}_M^+$ sub-homogeneous concave weight function.

Let $\chi\in\mathcal{W}:=\mathcal{W}^{-}\cup\mathcal{W}_M^{+}$. The set $\mathcal{E}_{\chi}(\Omega)$ is 
defined by 
$$\mathcal{E}_{\chi}(\Omega)=\{u\in\mpsh:\exists (u_j)_j\subset\mathcal{E}_0(\Omega),u_j\downarrow u~and~\sup_j \int_{\Omega}-\chi\circ u_j (\pp u_j)^n<+\infty\}.$$
For $u\in\mathcal{E}_{\chi}(\Omega)$, we let
$$E_{\chi}(u)=\int_{\Omega} -\chi(u)(\pp u)^n.$$
We note that the definition is well-defined. Indeed, we have 
\begin{proposition}\label{energy class in cegrell class}
    $$\mathcal{E}_{\chi}(\Omega)\subset\mathcal{E}(\Omega).$$
    In particular, the quaternionic Monge-Amp\`ere operator $(\pp u)^n$ is well-defined
    for any $u\in\mathcal{E}_{\chi}(\Omega)$ and $-\chi(u)\in L^1((\pp u)^n)$. Moreover, $(\pp u)^n$ puts no mass on quaternionic polar subsets.
\end{proposition}
\begin{proof}
     Let $r>0$ be a real number such that $\chi(-r)<0$.
    Then  we can choose an increasing
     function $\tilde{\chi}$ such that $\tilde{\chi}:\mathbb{R}^{-}\to\mathcal{R}^{-}$ is convex, $\tilde{\chi}\prime=\tilde{\chi}\prime\prime=0$ on $[-r,0]$ and $\tilde{\chi}\geq\chi$.
     We have that $\tilde{\chi}(v)\in\mpshn$. If $u\in\mathcal{E}_{\chi}(\Omega)$, then 
     by definition there exists a sequence $u_j\in\mathcal{E}_0(\Omega)$ decreasing to $u$ and satisfying
     $\sup_j\int_{\Omega}-\chi(u_j)(\pp u_j)^n<+\infty$. Then 
     \[
\sup_{j \in \mathbb{N}} \int_\Omega -\tilde{\chi}(u_1)(\pp u_j)^n \leq \sup_{j \in \mathbb{N}} \int_\Omega -\chi(u_j) (\pp u_j)^n < \infty.
\]
Since the function $\tilde{\chi} \circ u_1$ is a negative QPSH function, we claim that $u \in \mathcal{E}(\Omega)$. Indeed, let $G \Subset \Omega$ be a subdomain and consider the function
\[
 u_j^G := \sup \{ v \in \text{QPSH}^-(\Omega); v \leq u_j \text{ on } G \}.
\]
Here $u_j^G \in \mathcal{E}_0(\Omega)$ and $u_j^G \searrow u$ on $G$. Let $\psi \in \mathcal{E}_0(\Omega)$ be given such that $\tilde{\chi} \circ u_1 \leq \psi$. By integration by parts, we have
\[
\sup_{j \in \mathbb{N}} \int_\Omega -\psi (\pp u_j^G)^n \leq \sup_{j \in \mathbb{N}} \int_\Omega -\psi (\pp u_j)^n \leq \sup_{j \in \mathbb{N}} \int_\Omega -\tilde{\chi}(u_1) (\pp u_j)^n < \infty.
\]
Then
\[
\sup_{j \in \mathbb{N}} \int_\Omega (\pp u_j^G)^n \leq \left(- \sup_G \psi \right)^{-1} \sup_{j \in \mathbb{N}} \int_\Omega -\psi (\pp u_j^G)^n < \infty,
\]
which implies that $u \in \mathcal{E}(\Omega)$.

Therefore $(\pp u_j)^n$ weakly converges to $(\pp u)^n$ \cite[Theorem 3.2]{Wan20a}. It follows from
the upper semicontinuity of $u$ that $E_{\chi}(u)\leq\liminf_{j\to\infty}E_{\chi}(u_j)<\infty$. The 
last statement can be implied by Lemma \ref{100}.
\end{proof}
Moreover, we have the following theorem:
\begin{theorem}
    $$\mathcal{E}_{\chi}(\Omega)=\{u\in\mathcal{N}(\Omega)|\chi\circ u\in L^1((\pp u)^n)\}.$$
\end{theorem}
\begin{proof}
    The proof is a quaternionic analogue of that  of \cite[Theorem 2.7]{Ben11}.
\end{proof}
\begin{lemma}\label{100}
    If $u\in\mathcal{E}(\Omega)$, then for all Borel sets $B\subset\Omega\setminus\{u=-\infty\}$,
    $$\int_B(\pp u)^n=\lim_{j\to\infty}\int_{B\cap\{u>-j\}}(\pp u_j)^n,$$
    where $u_j:=\max(u,-j)$. In particular, the measure $(\pp u)^n$
    puts no mass on Q-polar sets $E\subset\{u>-\infty\}$.
\end{lemma}
\begin{proof}
    Since this result is in local nature, we can assume that $u\in\mathcal{F}(\Omega)$.
    For $s>0$ let $h_s:=\max(u/s+1,0)\in\mpsh$. Note that $h_s$ decrease to 
    the Borel function $\mathbf{1}_{\{u>-\infty\}}$ and $\{h_s=0\}=\{u\leq -s\}$. 
    By Proposition \ref{plurifine property} we have 
    $$h_s(\pp\max(u,-s))^n=h_s(\pp u)^n, \forall s>0.$$
    Observe that $h_s(\pp\max(u,-s))^n=h_s\mathbf{1}_{\{u>-s\}}(\pp u)^n$
    increase to $\mathbf{1}_{\{u>-\infty\}}\mu_u=\mathbf{1}_{\{u>-\infty\}}$ as $s\to +\infty$.   
    Similarly $h_s(\pp u)^n$ converges to $\mathbf{1}_{\{u>-\infty\}}(\pp u)^n$.
    Hence $\mu_u=\mathbf{1}_{\{u>-\infty}(\pp u)^n$, which ends the proof.

\end{proof}
We have the following comparison principle for functions in $\mathcal{E}_{\chi}(\Omega)$:
\begin{theorem}\label{cp for energy class}
    If $u\in\mathcal{E}(\Omega)$ and $v\in\mathcal{E}_{\chi}(\Omega)$ are such that 
    $(\pp v)^n\geq(\pp u)^n$, then $u\leq v$.
\end{theorem}
\begin{proof}
    One can follow the line in \cite[Theorem 4.5]{Ben09} combining the results in \cite{Wan20a}\cite{WZ15}.
\end{proof}
It is useful to introduce the following classes in terms of the speed of decrease of 
the capacity of sublevel sets.
\begin{definition}
    Let $\chi$ be a weight function defined as above. The class $\hat{\mathcal{E}_{\chi}}(\Omega)$
    is defined as
    $$\hat{\mathcal{E}_{\chi}}(\Omega):=\bigg\{u\in\mpshn:\bigg|\int_{0}^{+\infty}t^n\chi^{\prime}(-t) Cap(\{u<-t\})dt<+\infty\bigg\}.$$
\end{definition}

\begin{proposition}\label{energy class inclusion}
    We have
        \begin{itemize}
            \item[(1)] $\hat{\mathcal{E}}_{\chi}(\Omega)\subset\mathcal{E}_{\chi}(\Omega)$;
            \item[(2)] $\mathcal{E}_{\chi}(\Omega)\subset\hat{\mathcal{E}}_{\hat{\chi}}(\Omega)$, where $\hat{\chi}(t)=\chi(t/2)$;
            \item[(3)]  $ \mathcal{E}_{\chi}(\Omega)=\hat{\mathcal{E}}_{\chi}(\Omega)$, in particular, if $u\in\mathcal{E}_{\chi}$, then $(\pp u)^n$ puts no mass on Q-polar subsets;
        \item[(4)] The set $\mathcal{E}_{\chi}(\Omega)$ is a convex cone.
        \end{itemize}
\end{proposition}
\begin{proof}
    If $\chi\in\mathcal{W}_M^{+}$, one can adapt the approach in \cite{Sal23} to obtain the proof. As for $\chi\in\mathcal{W}^{-}$, see
    \cite[Proposition 2.4]{Ben11}. We only need to prove (4). However, this follows from (3).
\end{proof}
The following result corresponds to complex versions of \cite[Theorem 6]{Ben15} or \cite[Theorem 5.1]{Ben09} in the case of convex weight $\chi$ and 
\cite[Theorem 3.5]{Sal23} in the case of sub-homogeneous concave weight.
\begin{theorem}\label{characterization}
    Let $\mu$ be a positive Radon measure, and let $\chi \in \mathcal{W}$. The following conditions are equivalent:
    \begin{enumerate}
        \item there exists a \textit{unique} function $\varphi \in \mathcal{E}_{\chi}(\Omega)$ such that $\mu = (dd^c \varphi)^n$;
        \item\label{cond2} $\chi(\mathcal{E}_{\chi}(\Omega)) \subset L^1(d\mu)$;
        \item\label{cond3} there exists a constant $C > 0$ such that
        \[
        \int_{\Omega} -\chi \circ \psi \, d\mu \leq C,
        \]
        for all $\psi \in \mathcal{E}_0(\Omega)$, $E_{\chi}(\psi) \leq 1$;
        \item\label{cond4} there exists a positive constant $A$ such that
        \[
        \int_{\Omega} -\chi \circ \psi \, d\mu \leq A \max(1, E_{\chi}(\psi)), \quad \forall \psi \in \mathcal{E}_0(\Omega).
        \]
    \end{enumerate}

\end{theorem}
\begin{proof}
    The implication $(1)\implies (2)$ is obvious since $\mathcal{E}_{\chi}(\Omega)$ is convex and 
    $$\int_{\Omega}-\chi(u)(\pp v)^n\leq E_{\chi}(u+v)<+\infty, \forall u,v\in\mathcal{E}_{\chi}.$$
Now we prove $(2)\implies (3)$. By contradiction, let $(u_j)\in\mathcal{E}_0(\Omega)$
such that $E_{\chi}(u_j)\leq1$ and $$\int_{\Omega}-\chi(u_j)d\mu\geq 2^{3Mj}(M>1).$$
Let $$u:=\sum_j\frac{1}{2^{2j}}u_j. $$ 
Observe that $\{u<-s\}\cup_j\{u_j<-2^js\}$. Hence $$ Cap(\{u<-s\})\leq\sum_j Cap(\{u_j<-2^{j}s\}).$$
Now we have
\begin{equation*}
    \begin{aligned}
        \int_{0}^{\infty}s^n \chi\prime(-s) Cap(\{u<-s\})ds&\leq\int_{0}^{\infty}s^n\chi\prime(-s)\sum_{j}^{\infty} Cap(\{u_j<-2^{j}s\})ds\\
        &=\sum_{0}^{\infty}\frac{1}{2^{nj-n}}\int_{0}^{\infty}(2^{j-1})^n\chi\prime(-s)Cap(\{u_j\leq -2^js\})ds\\
        &\leq\sum_{0}^{\infty}\frac{1}{2^{nj-n}}\int_{0}^{\infty}\chi\prime(-s)\int_{\{u_j<-2^{j-1}s\}}(\pp u_j)^n ds\\
        &\leq\sum_{0}^{\infty}\frac{1}{2^{nj-n}}\int_{0}^{\infty}\chi\prime(-s)\int_{\{u_j<-s\}}(\pp u_j)^n ds\\
        &=\sum_{0}^{\infty}\frac{1}{2^{nj-n}}E_{\chi}(u_j)
        \leq\sum_{0}^{\infty}\frac{1}{2^{nj-n}}<\infty,
    \end{aligned}
\end{equation*}
where the second inequality uses Lemma \ref{cap estimate}. This proves that $u\in\hat{\mathcal{E}}_{\chi}(\Omega)$ and 
hence by Proposition \ref{energy class inclusion} we have $u\in\mathcal{E}_{\chi}(\Omega)$.
On the other hand, if $\chi$ is a convex weight, we have 
$$\int_{\Omega}-\chi(u)d\mu\geq\int_{\Omega}-\chi(\frac{1}{2^{2j}}u_j)\geq\frac{1}{2^{2j}}\int_{\Omega}-\chi(u_j)d\mu\geq 2^{Mj},\forall j\in\mathbb{N},$$
which yields a contradiction with (\ref{cond2}).
If $\chi$ is a sub-homogeneous concave weight, we have
 $$ (-\chi)(u_j)\leq 2^{2Mj}(-\chi)(2^{-2j}u_j)\leq 2^{2Mj}(-\chi)(u),$$
 where the first inequality follows from \ref{subh ineq} and the second follows from the definition of $u$.
 Therefore, $$\int_{\Omega}-\chi(u)d\mu\geq 2^{-2Mj}\int_{\Omega}-\chi(u_j)d\mu\geq 2^{Mj},$$
 which also yields a contradiction. 

 Now we move to prove $(3)\implies(4)$. Let 
 $\psi\in\mathcal{E}_0(\Omega)$. We may assume that $E_{\chi}(\psi)>1$, otherwise we are done. First we deal with the case that $\chi$ is convex.
 Define $\tilde{\psi}$ by $$\tilde{\psi}:=\frac{\psi}{E_{\chi}(\psi)^{1/n}}\in\mathcal{E}_0(\Omega).$$
 Indeed, since $\chi$ is increasing, we have 
 $$\int_{\Omega}-\chi(\tilde{\psi})(\pp \tilde{\psi})^n\leq\frac{1}{E_{\chi}(\psi)}\int_{\Omega}-\chi(\psi)(\pp \psi)^n=1.$$
 It follows from (\ref{cond3}) and the convexity of $\chi$ that
 $$\int_{\Omega}-\chi(\psi)d\mu\leq E_{\chi}(\psi)^{1/n}\int_{\Omega}(\frac{\psi}{E_{\chi}(\psi)^{1/n}})d\mu\leq CE_{\chi}(\psi)^{1/n}\leq CE_{\chi}(\psi).$$
 Hence we can choose $A$ in (\ref{cond4}) to be $C$ as above.
 
 If $\chi\in\mathcal{W}_M^{+}$, we use the method of QPSH envelop as in \cite{Sal23}. If 
 $E_{\chi}(\psi)\leq 2^{n+1}$, then $E_{\chi}(1/2\psi)\leq1$. If follows from \ref{subh ineq}
that $$\int_{\Omega}-\chi(\psi)d\mu\leq 2^M\int_{\Omega}-\chi(1/2\psi)d\mu\leq 2^M C.$$
So we may suppose that $E_{\chi}(\psi)\geq 2^{n+1}$ and denote 
by $\varepsilon = 1/(E_{\chi}(\psi))$. 
For $f = \chi^{-1}(\varepsilon \chi(\psi))$, consider the envelope
\[
P(f) = \sup \{ h \in \text{QPSH}(\Omega) : h \leq f \}.
\]
Obviously $P(f) \in \mathcal{E}_0(\Omega)$ since $f$ is upper-semicontinuous and satisfies $f \geq \psi$. 
Theorem \ref{mass conc} implies
\[
\mathbf{1}_{\{P(f) < f\}} (\pp P(f))^n = 0.
\]
It thus follows that
\[
E_{\chi}(P(f)) = \int_{\Omega} -\chi(P(f)) (\pp P(f))^n
\]
\[
= \int_{\{P(f)=f\}} -\chi(f) (\pp P(f))^n
\]
\[
= \varepsilon \int_{\Omega} -\chi(\psi) (\pp P(f))^n.
\]

By the proof of \cite[Theorem 3.4]{Sal23} for $\lambda = 1/2$, we get
\[
E_{\chi}(P(f))) \leq 2^n \varepsilon E_{\chi}(P(f)) + 2^{M+n}.
\]
This implies
\[
E_{\chi}(P(f)) \leq \frac{2^{M+n}}{1 - 2^n \varepsilon}\leq 2^{M+n+1},
\]
Therefore
\[
E_{\chi}(1/2^{M+1} P(f))=\int_{\Omega}-\chi(1/2^{M+1}P(f))(\pp \frac{P(f)}{2^{M+1}})^n<1
\]
by concavity of $\chi$.
Furthermore, since $P(f) \leq f = \chi^{-1}(\varepsilon \chi(\psi))$, we obtain
\[
\int_{\Omega} -\chi(\psi) d\mu \leq \varepsilon^{-1} \int_{\Omega} -\chi(P(f)) d\mu \leq 2^{M(M+1)} C E_{\chi}(\psi)
\]
by \ref{subh ineq}.
We conclude that, for all $\psi \in \mathcal{E}_0(\Omega)$, we have
\[
\int_{\Omega} -\chi(\psi) d\mu \leq 2^M C + 2^{M(M+1)} C E_{\chi}(\psi) \leq A \max(1, E_{\chi}(\psi)),
\]
where we choose $A= 2^{M(M+1)} + C$.

 Finally we prove $(4)\implies(1)$. Set ${\mu}=\frac{1}{2A}\mu$. Note that (\ref{cond4}) implies that 
 $\tilde{\mu}$ puts no mass on Quaternionic polar subsets by \cite[Theorem 1.1]{WK17}.
 Then it follows from \cite[Proposition 4.3]{Wan20a} that there exists $u\in\mathcal{E}_0(\Omega)$ and
 $f\in L_{loc}^1((\pp u)^n)$ such that ${\mu}=f(\pp u)^n$. Consider ${\mu}_j:=\min(f,j)(\pp u)^n$. By \cite[Lemma 4.7]{Wan20a}
 there exists a unique $\varphi_j\in\mathcal{E}_0(\Omega)$ such that $(\pp \varphi_j)^n={\mu}_j$. The comparison principle (e.g.\cite[Corollary 1.1]{WZ15})
 implies that $\varphi_j$ is a decreasing sequence. Set $\psi:=\lim_j\varphi_j$. We have
 $$\int_{\Omega}-\chi(\varphi_j)(\pp\varphi_j)^n\leq \frac{1}{2}\max(1,E_{\chi}(\varphi_j)),\forall j,$$
 which implies that $E_{\chi}(\varphi_j)\leq1$ is uniformly bounded with respect to $j$.
 
 On the other hand, either $\chi\in\mathcal{W}^{-}$ or $\chi\in\mathcal{W}_M^{+}$ gives that
 $\chi\prime(-2s)\leq M \chi\prime(-s)$ for some $M$ by definition. Combining it with Lemma \ref{cap estimate}, we have 
 $$\sup_j\int_{0}^{\infty}\chi\prime(-s)s^n Cap(\{\varphi_j<-s\})<+\infty.$$
Hence $\int_{0}^{\infty}\chi\prime(-s)s^n Cap(\{\varphi<-s\})<+\infty$ since 
$\varphi\leq\varphi_j$. Therefore we have $\varphi\in\mathcal{E}_{\chi}(\Omega)$.
We can conclude now by continuity of the complex Monge-Amp\`ere operator along
decreasing sequences (Proposition \ref{energy class in cegrell class}) that $(\pp \varphi)^n=\mu$.
The uniqueness of solution follows from the comparison principle for $\mathcal{E}_{\chi}$ (Theorem \ref{cp for energy class}).
 
\end{proof}

\end{document}